\documentclass{article}

\usepackage{amsmath}
\usepackage{amsthm}
\usepackage{amssymb}
\usepackage{graphicx}
\usepackage{natbib}
\setcitestyle{numbers}
\newtheorem{theorem}{Theorem}[section]
\newtheorem{definition}[theorem]{Definition}
\newtheorem{remark}[theorem]{Remark}
\numberwithin{subsection}{section}
\numberwithin{theorem}{subsection}
\numberwithin{equation}{subsection}

\title{Trigonometric Convexity for the Multidimensional Indicator after Ivanov}

\author{
\begin{minipage}{0.5\linewidth}
    Aleksandr Mkrtchyan$^{1,2}$\\
\small{1. Siberian Federal University, Russia\\
    2. Institute of Mathematics NAS, Armenia\\\\}
\end{minipage}
\hfill
\begin{minipage}{0.5\linewidth}
    Armen Vagharshakyan$^{3,4}$\\
\small{    3. Institute of Mathematics NAS, Armenia \\
4. Regional Mathematical Center of the\\
Southern Federal University, Russia}
\\
\end{minipage}
}
\date{}
\begin{document}
\maketitle
\abstract{ Multidimensional indicator after Ivanov is a generalization of the notion of indicator, that is well-known for analytic functions in one complex variable, to analytic functions in several complex variables.
   We prove an analogue of trigonometric convexity for it. Additionally, we show that our estimate is sharp.
    The proof is based on the multidimensional analogue of the sectorial Fourier inversion formula.
}

\section{ Introduction.}
\subsection{ Main results.}
We review the notion of indicator for analytic functions in one complex variable (see section \ref{historyone}) and its multidimensional analogues for functions in several complex variables (see section \ref{historymany}); in particular, 
we treat
the notion of multidimensional indicator after Ivanov. 

As explained in remark \ref{convexity}, our main theorem \ref{ivanov} establishes an analogue of trigonometric convexity for the multidimensional indicator after Ivanov.
We preface formulation of  theorem \ref{ivanov} by the following definitions:

\noindent
\begin{definition}\label{opensector}
	Denote by $\Delta_{\alpha_j}\subset \mathbb{C}$ the open sector determined by the angle $ 0<\alpha_j<\pi/2$ as follows:
	\begin{equation*}
	\Delta_{\alpha_j}=\lbrace z_j\in\mathbb{C}\setminus\lbrace 0\rbrace\colon \left|\arg\left(z_j\right)\right|<\alpha_j\rbrace.
	\end{equation*}
\end{definition}
\begin{definition}
	Recall that a function $f$ is of  finite exponential type $\left(h_1,\dots, h_n\right)$ in $\Delta_{\alpha_1}\times\dots\times \Delta_{\alpha_n}$ if
	for any $\varepsilon>0$ there exists a constant $k_{\varepsilon}\geq 0$ such that
	\begin{equation}\label{h}
		\left|f\left(z_1,\dots,z_n\right)\right|\leq k_{\varepsilon}e^{(h_1+\varepsilon)\left|z_1\right|+\dots+\left(h_n+\varepsilon\right)\left|z_n\right|},\quad \text{ for all } z_j\in \Delta_j,\; 1\leq j\leq n.
	\end{equation}
Note that in this article we tacitly assume that $h_1,\dots,h_n\geq 0$. 
\end{definition}
\begin{definition}
	Denote by $Exp\left(\alpha_1,\dots,\alpha_n\right)$ the class of functions $f$ that are analytic and 
	are
	of finite exponential type in $\Delta_1\times \dots\times\Delta_n.$
\end{definition}
\begin{definition}
 Following Ivanov 
  \cite{I}, for $f \in Exp\left(\alpha_1,\dots,\alpha_n\right)$ denote by  $T_{f}\left(\vec \theta\right)$ 
  the following:
\begin{align}\label{T}
T_{f}\left(\vec \theta\right)=
\{
&
\vec \nu\in \mathbb R^n: 
\ln{\left|f\left(\vec re^{i\vec \theta}\right)\right|}\leq \nu_1 r_1+...+\nu_nr_n+C_{\vec \nu,\vec \theta},
\\\nonumber
&
\text{ for some } C_{\vec \nu,\vec\theta}, \text{ for all }\vec r\in \mathbb R^n_+
\}
,
\end{align}
here $\vec re^{i\vec \theta}$ 
denotes
the vector $\left(r_1e^{i\theta_1},...,r_ne^{i\theta_n}\right).$
\end{definition}
\begin{theorem}\label{ivanov}
 Let a function $f \in Exp\left(\alpha_1,\alpha_2\right)$ and  
  the numbers $A^+_1,A^+_2,A^-_1,A^-_2$ satisfy:
 \begin{align*}
 \left(A^+_1,A^+_2\right) \in \overline T_f\left(\alpha_1,\alpha_2\right), \quad
 \left(A^-_1,A^-_2\right) \in \overline T_f\left(-\alpha_1,-\alpha_2\right)
 ,
 \\
 \left(A^+_1,A^-_2\right) \in \overline T_f\left(\alpha_1,-\alpha_2\right), \quad
  \left(A^-_1,A^+_2\right) \in \overline T_f\left(-\alpha_1,\alpha_2\right).
 \end{align*}
 Then
 we have
 \begin{equation*}
 \left(C_1,C_2\right)\in \overline T_f\left(\theta_1,\theta_2\right),
 \end{equation*}
 where  the constants $C_1,C_2$ 
 are determined
 by
  the following formulas:
 \begin{align}\label{c}
 C_1\sin\left(2\alpha_1\right)=A^+_1\sin\left(\theta_1+\alpha_1\right)+A^-_1\sin\left(\alpha_1-\theta_1\right),
 \\
 C_2\sin\left(2\alpha_2\right)=A^+_2\sin\left(\theta_2+\alpha_2\right)+A^-_2\sin\left(\alpha_2-\theta_2\right).
  \end{align}
\end{theorem}

Theorem \ref{ivanov} can be paraphrased as follows:
\begin{remark}\label{main}
	 Let a function $f \in Exp\left(\alpha_1,\alpha_2\right)$ and  
	the numbers $A^+_1,A^+_2,A^-_1,A^-_2$ satisfy:
\begin{align}
	\nonumber
	&
	\left|f\left(z_1,z_2\right)\right|\leq k_{\varepsilon}e^{(A^+_1+\varepsilon)\left|z_1\right|+\left(A^+_2+\varepsilon\right)\left|z_2\right|},\quad \text{for }\arg(z_1)=\alpha_1,\;\arg\left(z_2\right)=\alpha_2,
	\\\nonumber
	&
	\left|f\left(z_1,z_2\right)\right|\leq k_{\varepsilon}e^{(A^+_1+\varepsilon)\left|z_1\right|+\left(A^-_2+\varepsilon\right)\left|z_2\right|},\quad \text{for }\arg(z_1)=\alpha_1,\;\arg\left(z_2\right)=-\alpha_2,
	\\\nonumber
	&
	\left|f\left(z_1,z_2\right)\right|\leq k_{\varepsilon}e^{(A^-_1+\varepsilon)\left|z_1\right|+\left(A^+_2+\varepsilon\right)\left|z_2\right|},\quad \text{for }\arg(z_1)=-\alpha_1,\;\arg\left(z_2\right)=\alpha_2,
	\\\label{ab}
	&
	\left|f\left(z_1,z_2\right)\right|\leq k_{\varepsilon}e^{(A^-_1+\varepsilon)\left|z_1\right|+\left(A^-_2+\varepsilon\right)\left|z_2\right|},\quad \text{for }\arg(z_1)=-\alpha_1,\;\arg\left(z_2\right)=-\alpha_2,
\end{align}
where by the value of the function $f$ on 
the mentioned
rays we mean $f$'s non-tangential limit (see \cite{V}, remark 2.1 for details).
Then
 for any $\varepsilon>0$ there exists a constant $k^{\prime}_{\varepsilon}>0$ such that
\begin{equation}\label{f_estimate}
	\left|f\left(z_1,z_2\right)\right|\leq
	k^{\prime}_{\varepsilon}e^{(C_1+\varepsilon)\left|z_1\right|+\left(C_2+\varepsilon\right)\left|z_2\right|},\quad \text{for }
	\arg\left(z_1\right)=\theta_1,\;\arg\left(z_2\right)=\theta_2,
\end{equation}
where 
  the constants $C_1,C_2$ are determined by the following formulas:
\begin{align*}
C_1\sin\left(2\alpha_1\right)=A^+_1\sin\left(\theta_1+\alpha_1\right)+A^-_1\sin\left(\alpha_1-\theta_1\right),
\\
C_2\sin\left(2\alpha_2\right)=A^+_2\sin\left(\theta_2+\alpha_2\right)+A^-_2\sin\left(\alpha_2-\theta_2\right).
\end{align*}
\end{remark}
In section \ref{ivanov_proof} we derive theorem \ref{ivanov} from 
theorem \ref{fourier_inversion}
 that serves as a two-dimensional analogue of the Fourier inversion formula (see \cite{V}, theorem 1.2 and \cite{DzA}) for functions of exponential type in a sector.
We preface formulation of  theorem \ref{fourier_inversion} by the following definitions:
\begin{definition}\label{laplastrans}
Let	$f \in Exp\left(\alpha_1,\alpha_2\right)$ and   
let 
the numbers $A^+_1,A^+_2,A^-_1,A^-_2$ satisfy \ref{ab}.  Define the function $m$  as $f$'s two-dimensional concatenated Laplace transform. Namely, the domain of  function $m$ is  the Cartesian product
$\Omega_1\times \Omega_2,$ 
where
\begin{align*}&
\Omega_1=\Omega_1^+\cup\Omega_1^-,
\\\nonumber &
\Omega_2=\Omega_2^+\cup\Omega_2^-,
\end{align*}
and in turn 
\begin{align}\label{omega_pm}
&
\Omega_1^+=\left\lbrace \omega_1 \colon Re\left(\omega_1 e^{i\alpha_1}\right)<-A^+_1\right\rbrace, 
\\\nonumber &
\Omega_1^-=\left\lbrace \omega_1 \colon Re\left(\omega_1 e^{-i\alpha_1}\right)<-A^-_1\right\rbrace, 
\\\nonumber &
\Omega_2^+=\left\lbrace \omega_2 \colon Re\left(\omega_2 e^{i\alpha_2}\right)<-A^+_2\right\rbrace, 
\\\nonumber &
\Omega_2^-=\left\lbrace \omega_2 \colon Re\left(\omega_2 e^{-i\alpha_2}\right)<-A^-_2\right\rbrace, \end{align}
\begin{figure}\label{fig_omega}
	\centering
	\includegraphics{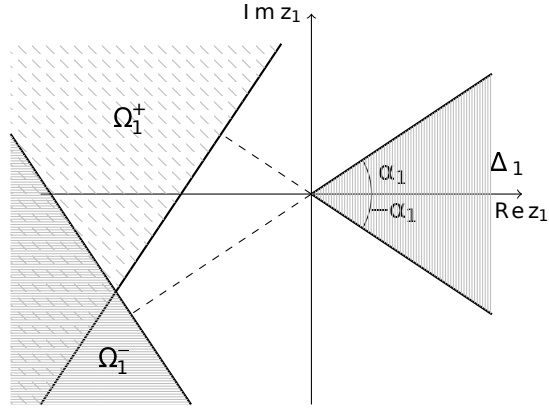}
	\caption{The set $\Omega$}
\end{figure}
 (see figure \ref{fig_omega} on page \pageref{fig_omega}). The function $m$ is defined
on  $\Omega_1\times \Omega_2$ 
  by  the following four formulas:
\begin{align}\label{m}
&
m\left(\omega_1,\omega_2\right)=\frac{-1}{4\pi^2}
\int_{e^{i\alpha_2}[0,+\infty)}\int_{e^{i\alpha_1}[0,+\infty)}
f\left(\zeta_1,\zeta_2\right)
e^{\omega_1\zeta_1+\omega_2\zeta_2}
d\zeta_1d\zeta_2
,\\\nonumber &\quad\quad\quad\quad 
\text{for }
\omega_1\in\Omega_1^+
,\;
\omega_2\in\Omega_2^+,
\\\nonumber
&
m\left(\omega_1,\omega_2\right)=\frac{-1}{4\pi^2}
\int_{e^{i\alpha_2}[0,+\infty)}\int_{e^{-i\alpha_1}[0,+\infty)}
f\left(\zeta_1,\zeta_2\right)
e^{\omega_1\zeta_1+\omega_2\zeta_2}
d\zeta_1d\zeta_2
,\\\nonumber & \quad\quad\quad\quad\text{for }
\omega_1\in\Omega_1^+
,\;
\omega_2\in\Omega_2^-,
\\\nonumber
&
m\left(\omega_1,\omega_2\right)=\frac{-1}{4\pi^2}
\int_{e^{-i\alpha_2}[0,+\infty)}\int_{e^{i\alpha_1}[0,+\infty)}
f\left(\zeta_1,\zeta_2\right)
e^{\omega_1\zeta_1+\omega_2\zeta_2}
d\zeta_1d\zeta_2
,\\\nonumber &\quad\quad\quad\quad \text{for }
\omega_1\in\Omega_1^-
,\;
\omega_2\in\Omega_2^+,
\\\nonumber
&
m\left(\omega_1,\omega_2\right)=\frac{-1}{4\pi^2}
\int_{e^{-i\alpha_2}[0,+\infty)}\int_{e^{-i\alpha_1}[0,+\infty)}
f\left(\zeta_1,\zeta_2\right)
e^{\omega_1\zeta_1+\omega_2\zeta_2}
d\zeta_1d\zeta_2
,\\\nonumber & \quad\quad\quad\quad\text{for }
\omega_1\in\Omega_1^-
,\;
\omega_2\in\Omega_2^-.
\end{align}
	\end{definition}
\begin{definition}\label{definitiongamma}
For a class of functions $Exp\left(\alpha_1,\alpha_2\right),$ denote by
$\Gamma_1$ and $\Gamma_2$  the curves given by the following parametrizations:
\begin{align}\label{gamma_1}
&
\gamma_1\colon\mathbb{R}\rightarrow\mathbb{C},\nonumber
\\
&
\gamma_1(t)=p_1-ie^{i\alpha_1}|t|,\quad\text{for }t\in (-\infty,0],\nonumber
\\
&
\gamma_1(t)=p_1+ie^{-i\alpha_1}|t|,\quad\text{for }t\in (0,+\infty],
\end{align}
(see figure \ref{fig_gamma} on page \pageref{fig_gamma}), and the real number $p_1$ appearing in parametrization \eqref{gamma_1}, is chosen in such a way that it satisfies inequality
\begin{equation}\label{p_1}
p_1\cos(\alpha_1)<-h_1
\end{equation}
\begin{figure}\label{fig_gamma}
	\centering
	\includegraphics{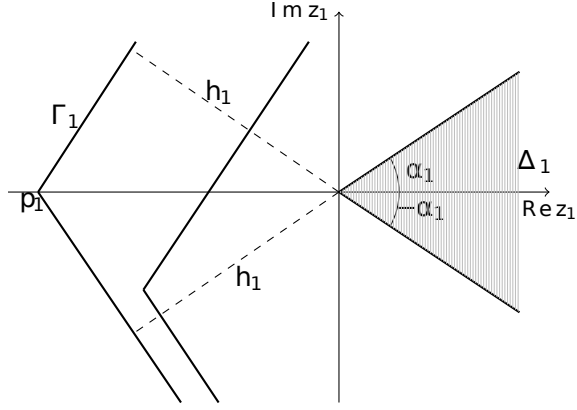}
	\caption{The curve $\Gamma_1$}
\end{figure}
and correspondingly
the curve $\Gamma_2$ is parameterized by
\begin{align}\label{gamma_2}
&
\gamma_2\colon\mathbb{R}\rightarrow\mathbb{C},\nonumber
\\
&
\gamma_2(t)=p_2-ie^{i\alpha_2}|t|,\quad\text{for }t\in (-\infty,0],\nonumber
\\
&
\gamma_2(t)=p_2+ie^{-i\alpha_2}|t|,\quad\text{for }t\in (0,+\infty],
\end{align}
and the real number $p_2$ appearing in parametrization \eqref{gamma_1}, is chosen in such a way that it satisfies inequality
\begin{equation}\label{p_2}
p_2\cos(\alpha_2)<-h_2,
\end{equation}
\end{definition}
\begin{theorem}\label{fourier_inversion}
For a function $f \in Exp\left(\alpha_1,\alpha_2\right) $  the following Fourier inversion formula  holds: 
\begin{equation}\label{ivanov_formula}
    f\left(z_1,z_2\right)=\int_{\Gamma_1}\int_{\Gamma_2}m\left(\omega_1,\omega_2\right)
    e^{-\omega_1 z_1-\omega_2z_2}
    d\omega_2d\omega,\quad \text{for }z_1\in \Delta_1,\;z_2\in \Delta_2.
\end{equation}
\end{theorem}
Theorem \ref{fourier_inversion} is proved in section \ref{fourier_inversion_proof}.

We complement theorem \ref{ivanov} by the following remarks:
\begin{remark}
Theorem \ref{ivanov} is stated for functions of two complex variables. The corresponding theorem for functions of $n$ complex variables also holds.
\end{remark}
\begin{remark}\label{sharp}
Theorem \ref{ivanov} is sharp: that is, there exists a function $f$ for whom: the assumptions of theorem \ref{ivanov} are satisfied, and the inequality \eqref{f_estimate} is an equality.
\end{remark}
Remark \ref{sharp} is proved in section \ref{sharp_section}.
\begin{remark}
Integral representation \eqref{ivanov_formula} in \ref{fourier_inversion} could be compared to multivariate integrals of Mellin-Barnes type \cite{FG},\cite{PTY}.
\end{remark}
\subsection{ Indicators of functions in one complex variable.}\label{historyone}
We say that an entire function $f\colon\mathbb{C}\rightarrow\mathbb{C}$ is of exponential type if
\begin{equation}\label{indicatorone}
\limsup_{z \to \infty}  \frac{\ln^+|f(z)|}{|z|}<+\infty.
\end{equation}
The notion of indicator of an entire function of exponential type was introduced by Phragmen and Lindeleph \cite{PhL},\cite{Bi} as follows:
\begin{equation}\label{indicator}
h_{f}(\alpha):=\limsup_{r \to \infty}  \frac{\ln\left|f\left(re^{i\alpha}\right)\right|}{r}, \;\; \alpha \in \mathbb{R}.
\end{equation}
The indicator $h_f(\alpha)$ describes the growth of function  $f$ along the ray $e^{i\alpha}[0,+\infty).$
It follows from definition \eqref{indicator} that the indicator 
 $h_{f}(\alpha)$ is a real-valued $2\pi$-periodic function.
It also follows that the indicator of the product of two functions does not exceed the sum of the indicators of the factors, 
\begin{equation*}
h_{fg}(\alpha)\leq h_{f}(\alpha)+h_{g}(\alpha),
\end{equation*}
and that the indicator of the sum of two functions does not exceed the larger of the two indicators,
\begin{equation*}
h_{f+g}(\alpha)\leq \max\left(h_{f}(\alpha),h_{g}(\alpha)\right).
\end{equation*}
 One of the main properties of the indicator
$h_{\theta}(\alpha)$ is its \textit{trigonometric convexity} \cite{Lev},\cite{Bo}:
if $\alpha_1<\alpha<\alpha_2$ and  $\alpha_2-\alpha_1<\pi,$ then the following inequality holds:
\begin{equation}\label{trig_convex}
h_{f}(\alpha)\sin{(\alpha_2-\alpha_1)}\leq h_{f}(\alpha_1)\sin{(\alpha_2-\alpha)}+h_{f}(\alpha_2)\sin{(\alpha-\alpha_1)}.
\end{equation}
The following property follows from the trigonometric convexity \eqref{trig_convex} \cite{Lev}: 
if the indicator is bounded on an open interval, then it is continuous.
The latter claim does not hold for a closed interval.

The notion of indicator is known to be an important tool in some methods regarding to finding analytic continuation
\cite{P}, \cite{C},\cite{C1},\cite{Ar85}. 
In particular, the indicator appears in problems relating to analytic continuation of power series via interpolation of coefficients. It also plays a role in problems relating to localization of singularities of power series
 \cite{ALM}.
 
  Let 
\begin{equation*}
f(z)=\sum_{k=0}^\infty {\frac{a_k}{k!}z^k},
\end{equation*}
be the power-series representation of the entire function
 $f.$ 
 Consider its Borel transform defined by
the Laurent series in the following way:
\begin{equation}\label{borel}
g(\omega)=\sum_{k=0}^\infty {a_k \omega^{-k-1}}.
\end{equation}

The interrelation between the set of singularities of
 $g$ and the indicator of   $f$ is described by Polya's theorem  \cite{Leo}, \cite{Ro}:
 \begin{theorem}[Polya]\label{P}
Let $f$ be an entire function of exponential type.
Denote by $K\subset \mathbb{C}$ the convex set whose support function 
\begin{equation*}
k(\theta)=\sup_{\omega\in K}Re\left(\omega e^{-i\theta}\right)
\end{equation*}
is determined by $f$'s indicator as follows:
\begin{equation}
    k(-\theta)=h_f(\theta).
\end{equation}
Then $f$ can be restored by
\begin{equation}
    f(z)=\frac{1}{2\pi i}\int_{\Gamma}g(\omega)e^{z\omega}d\omega,
\end{equation}
where $\Gamma$ is a closed contour containing the set $K,$ and
$g$ \eqref{borel} is the Borel transform of $f.$ Additionally, $K$ is the smallest convex set such that $g$ is analytic in $\mathbb{C}\setminus K.$
\end{theorem}
\subsection{ Indicators of functions in several complex variables.}\label{historymany}
The works of Ronkin \cite{Ro}, Lelon \cite{LG}, Levin \cite{Lev}, Ivanov \cite{I}, Kiselman \cite{Ki} and others  are dedicated to exploring multidimensional analogues of the indicator function.   
We say that  $f\left(\vec z\right)=f(z_1,...,z_n)$ is an $n$-valued entire function of exponential type if $f$ is holomorphic in
$\mathbb{C}^n,$ 
\begin{equation*}
    f\left(\vec z\right) \in \mathcal O (\mathbb C^n),
\end{equation*}
and if there exist constants  $k, \sigma_1,...,\sigma_n$ such that
\begin{equation}
    \left|f\left(\vec z\right)\right|\leq ke^{\sigma_1|z_1|+...+\sigma_n|z_n|},\quad \text{for }\vec z\in\mathbb{C}^n.
\end{equation}
In the case of several complex variables, different characteristics of an analytic function's growth in directions have been introduced. 
For example, introduce
the \textit{radial indicator } of the function $f$ as follows \cite{Ro}, \cite{LG}: $$L_r(z,f)=\limsup_{t \to \infty}\frac{\ln|f(tz)|}{t}, \;\; z\in\mathbb{C}^n$$
and correspondingly introduce regularization of the function $L(z,f)$ as follows:
\begin{equation}
L^{\ast}(z,f)=\limsup_{z' \to z}L(z,f)
\end{equation}
 We will refer to $L^{\ast}(z,f)$ as \textit{regularized radial indicator} of the function $f$.
Just as in the on-dimensional case the function $L^{\ast}(z,f)$ is semi-continuous from above. Note that the regularized radial indicator is a plurisubharmonic function in $\mathbb{C}^n.$
Consequently, more information about its properties may be found in works related to plurisubharmonic functions in potential theory \cite{Ki} and their different generalizations \cite{Ko}, \cite{MaS}.

Another charateristic of an analytic functions growth was introduced by Ivanov; namely Ivanov \cite{I} has introduced the set $T_{f}\left(\vec \theta\right)$ as follows: 
\begin{align}\label{Ts}
T_{f}\left(\vec \theta\right)=
\{&
\vec \nu\in \mathbb R^n: 
\ln{\left|f\left(\vec re^{i\vec \theta}\right)\right|}\leq \nu_1 r_1+...+\nu_nr_n+C_{\vec \nu,\vec \theta},
\\\nonumber
&
\text{ for some } C_{\vec \nu,\vec\theta}, \text{ for all }\vec r\in \mathbb R^n_+
\},
\end{align}
here $\vec re^{i\vec \theta}$ is the vector $(r_1e^{i\theta_1},...,r_ne^{i\theta_n}).$
The set $T_{f}\left(\vec \theta\right)$ implicitly reflects the notion of an indicator of an entire function.

		 For example, for the closure of the set $T_\varphi\left(\theta_1,\theta_2\right)$ defined for the  function $\varphi(\zeta_1,\zeta_2)=\cos(\zeta_1\zeta_2)^{\frac{1}{2}}$  of exponential type, 
we have the following \cite{Mk}:
		\begin{equation*}
		\overline{T}_{\varphi}\left(\theta_1,\theta_2\right)=\left\lbrace
		\vec \nu\in\mathbb R^2: \nu_1\nu_2=\frac{1}{4}\left|\sin\left(\frac{\theta_1+\theta_2}{2}\right)\right|^2,\;\;\nu_1\geq0, \;\; \nu_2\geq0 
		\right\rbrace.
		\end{equation*}
There exist many more analogues of indicator for functions in several complex variables. However, for none of them a property resembling trigonometric convexity was obtained. 
\begin{remark}\label{convexity}
Following our main result \ref{ivanov}, 
we now formulate trigonometric convexity for multidimensional indicator after Ivanov.

 Let a function $f \in Exp\left(\alpha_1,\dots,\alpha_n\right)$ and  
the numbers $A^+_1,A^-_1 \dots,A^+_n,A^-_n$ satisfy
\begin{align*}
\left(A^{l_1}_1, \dots,A^{l_n}_n\right) \in \overline T_f\left({l_1}\alpha_1, \dots,{l_n}\alpha_n\right), 
\end{align*}
where $l_j=\pm, \;$ $\; j=1,\dots,n.$  
Then
\begin{equation*}
\left(C_1,\dots,C_n\right)\in \overline T_f\left(\theta_1,\dots,\theta_n\right),
\end{equation*}
where  the constants $C_1,\dots,C_n$ determine from the following formulas:
\begin{align*}
C_j\sin\left(2\alpha_j\right)=A^+_j\sin\left(\theta_j+\alpha_j\right)+A^-_j\sin\left(\alpha_j-\theta_j\right), \;\; j=1,...,n.
\end{align*}
\end{remark}
\section{ Two-dimensional sectorial Fourier inversion formula.}\label{fourier_inversion_proof}
\subsection{ Two-dimensional concatenated Laplace transform}
\noindent
Due to \eqref{omega_pm} for any pair of complex numbers $\left(\omega_1,\omega_2\right)\subset \Omega_1^+\times\Omega_2^+$ we can
pick $\varepsilon>0$ so that the following inequalities are satisfied:
\begin{align}\label{close}
&
Re\left(\omega_1 e^{i\alpha_1}\right)+\varepsilon<-A^+_1,
\\\nonumber &
Re\left(\omega_2 e^{i\alpha_2}\right)+\varepsilon<-A^+_2.
\end{align}
Then we have the following estimate on the function $m$ defined by the first formula in \eqref{m}:
\begin{align}\label{mestimate}
&
   \left|m\left(\omega_1,\omega_2\right)\right|
  \overset{\eqref{m}}{\leq}
 \frac{1}{4\pi^2}
\int_{e^{i\alpha_1}[0,+\infty)}
\int_{e^{i\alpha_2}[0,+\infty)}
\left|f\left(\zeta_1,\zeta_2\right)\right|
\left|e^{\omega_1\zeta_1+\omega_2\zeta_2}\right|
\left|d\zeta_2\right|\left|d\zeta_1\right|
\leq
\\\nonumber &
\overset{\ref{ab}}{\leq}
     \frac{k_{\varepsilon}}{4\pi^2}
  \left[
\int_{e^{i\alpha_1}[0,+\infty)}
e^{
\left(
A^+_1+\varepsilon+Re\left(\omega_1 e^{i\alpha_1}\right)
\right)
\left|\zeta_1\right|
}
\left|d\zeta_1\right|
\right]
\cdot
\\\nonumber &
\cdot
\left[
\int_{e^{i\alpha_2}[0,+\infty)}
e^{
\left(
A^+_2+\varepsilon+Re\left(\omega_2e^{i\alpha_2}\right)
\right)\left|\zeta_2\right|
}
\left|d\zeta_2\right|\right]
=
\\\nonumber &
\overset{\eqref{close}}{=}
 \frac{k_{\varepsilon}}{4\pi^2}\cdot
 \frac{-1}{A^+_1+\varepsilon+Re\left(\omega e^{i\alpha_1}\right)}\cdot\frac{-1}{A^+_2+\varepsilon+Re\left(\omega_2 e^{i\alpha_2}\right)}<+\infty.
\end{align}
Due to estimate \eqref{mestimate} the first of four double integrals in \eqref{m} is absolutely convergent. And that double integral determines a function that is analytic in two complex variables on the set $\Omega_1^+\times\Omega_2^+.$ The same claims are true for the other three double integrals in the definition \eqref{m}.
We now address the formal ambiguity in definition \eqref{m} of function $m,$ arising from the fact that some of the four domains that appear in \eqref{m} might  intersect. Due to lemma 4.1 in \cite{V} (a lemma that: uses the exponential estimate \ref{h} on the function $f,$  and that is based on application of the Phragmen-Lindeloef maximum principle), the first and the second double integrals in \eqref{m} are equal on the intersection of their corresponding domains: $\Omega_1^+\times \Omega_2^+$ and $\Omega_1^+\times \Omega_2^-.$
Due to \eqref{mestimate} Fubini's theorem applies to each of the four double integrals in \eqref{m}. By changing the order of integration in the first and third double integrals and referring lemma 4.1 in\cite{V},
the first and the third double integrals are equal on the intersection of their corresponding domains. As for the first and the fourth double integrals, the intersection of their corresponding domains lies within the intersection of any two of the four domains. Hence in there, the first integral equals the second integral, and in turn, the second integral equals the fourth integral. Consequently, any two of the four definitions of the function $m$ in \eqref{m} are equivalent on the intersection of their corresponding domains.
The function $m$ defined by \eqref{m} is analytic on $\Omega_1\times\Omega_2$ in each of its variables.
\begin{remark}\label{boundedaway}
Thanks to the estimate \eqref{mestimate}, the function $m$ is bounded on any subset of $\Omega_1^+\times\Omega_2^+$ that is bounded away from its boundary $\partial\left(\Omega_1^+\times\Omega_2^+\right).$ Similarly, one can prove that
the function $m$ is bounded on any subset of $\Omega_1\times\Omega_2$ that is bounded away from the boundary $\partial\left(\Omega_1\times\Omega_2\right).$
\end{remark}
 \subsection{ Verification of the sectorial Fourier inversion formula.}
Due to the condition \ref{h},  for any $z_2\in \Delta_2$ the Fourier inversion formula (see \cite{V}, theorem 1.2)  applies to the function $z_1\rightarrow f\left(z_1,z_2\right):$ that is, we have
 \begin{equation}\label{fourier_first}
     f(z_1,z_2)=\int_{\Gamma_1}g\left(\omega_1,z_2\right)e^{-\omega_1 z_1}d\omega_1
     ,\quad\text{for } z_2\in \Delta_2,
 \end{equation}
where the curve $\Gamma_1$ is defined by \eqref{gamma_1}, and the function $\omega_1\rightarrow g\left(\omega_1,z_2\right)$ is well-defined by the following formulas:
\begin{align}\label{g}
&
    g\left(\omega_1,z_2\right)=\frac{1}{2\pi i}
\int_{e^{-i\alpha_1}[0,+\infty)}
f\left(\zeta_1,z_2\right)
e^{\omega_1 \zeta_1}
d\zeta_1
,\quad \text{for } Re\left(\omega_1 e^{-i\alpha_1}\right)<-h_1,
\\\nonumber
&
g\left(\omega_1,z_2\right)=\frac{1}{2\pi i}
\int_{e^{i\alpha_1}[0,+\infty)}
f\left(\zeta_1,z_2\right)
e^{\omega_1 \zeta_1}
d\zeta_1
,\quad \text{for } Re\left(\omega_1 e^{i\alpha_1}\right)<-h_1.
\end{align}
We now derive the estimate \eqref{g_uniform} for function $g\left(\omega_1,z_2\right),$ that is uniform in its first variable $\omega_1$ for $\omega_1\in\Gamma_1$
and that depends exponentially on its second variable $z_2.$ 

Due to inequality \eqref{p_1} imposed on the real number $p_1,$ we can pick $\varepsilon>0$ such that
\begin{equation}\label{epsilon}
    h_1+\varepsilon+p_1\cos\left(\alpha_1\right)<0.
\end{equation}
Let $\omega_1\in\Gamma_1.$ Then due to parameterization \eqref{gamma_1} we have $\omega_1=\gamma_1(t),\; t\in\mathbb{R}.$ Specifically, assume that
$
    t\geq 0,
$
so that the second formula of the two formulas \eqref{g} holds for $g\left(\gamma_1(t),z_2\right).$ 
We estimate
\begin{align}\label{g_uniform}
&
\left|g\left(\gamma_1(t),z_2\right)\right|
\overset{\eqref{g}}{\leq}
\frac{1}{2\pi}
\int_{e^{i\alpha_1}[0,+\infty)}
\left|f\left(\zeta_1,z_2\right)\right|
\cdot
\left|e^{\gamma_1(t) \zeta_1}\right|
\left|d\zeta_1\right|\leq
\\\nonumber
&
\overset{\ref{h}}{\leq}
\frac{1}{2\pi}
\int_{e^{i\alpha_1}[0,+\infty)}
k_{\varepsilon}e^{\left(h_1+\varepsilon\right)\left|\zeta_1\right|
+\left(h_2+\varepsilon\right)\left|z_2\right|}
\cdot
e^{Re\left(\gamma_1(t) \zeta_1\right)}
\left|d\zeta_1\right|=
\\\nonumber 
&
=
\frac{k_{\varepsilon}e^{\left(h_2+\varepsilon\right)\left|z_2\right|}}{2\pi}
\int_{e^{i\alpha_1}[0,+\infty)}
e^{\left(h_1+\varepsilon+Re\left(\gamma_1(t) e^{i\alpha_1}\right)\right)\left|\zeta_1\right|}
\left|d\zeta_1\right|=
\\\nonumber
&
\overset{\eqref{gamma_1}}{=}
\frac{k_{\varepsilon}e^{\left(h_2+\varepsilon\right)\left|z_2\right|}}{2\pi}
\int_{e^{i\alpha_1}[0,+\infty)}
e^{\left(h_1+\varepsilon+p_1\cos\left(\alpha_1\right)\right)\left|\zeta_1\right|}
\left|d\zeta_1\right|=
\\\nonumber
&
\overset{\eqref{epsilon}}{=}
\left[
\frac{k_{\varepsilon}}{2\pi}
\cdot
\frac{-1}{h_1+\varepsilon+p_1\cos\left(\alpha_1\right)}
\right]
\cdot
e^{\left(h_2+\varepsilon\right)\left|z_2\right|}.
\end{align}
And we would get the same estimate on $g\left(\gamma_1(t),z_2\right)$ if we assumed that $t<0.$
Due to estimate \eqref{g_uniform},  for any $-\infty<t<+\infty,$ the Fourier inversion formula (see \cite{V}, theorem 1.2) applies to the function $z_2\rightarrow g\left(\gamma_1(t),z_2\right),$ and we have
 \begin{equation}\label{second}
    g\left(\gamma_1(t),z_2\right)=\int_{\Gamma_2}
    m\left(\gamma_1(t),\omega_2\right)e^{-\omega_2 z_2}d\omega_2
     ,\quad\text{for } z_2\in \Delta_2,
 \end{equation}
 where the curve $\Gamma_2$ is parameterized by \eqref{gamma_2},
and for any $\infty <t<+\infty,$ the function $\omega_2\rightarrow m\left(\gamma_1(t),\omega_2\right)$ is well-defined by the formulas \eqref{m}.
\begin{remark}
Note that due to the uniform bound \eqref{g_uniform}, the curve $\Gamma_2,$ involved in formula \eqref{second}, is defined by formula \eqref{gamma_2} in such a way that it does not depend on the choice of $\omega_1=\gamma_1(t).$
\end{remark}
By combining the Fourier inversion formulas \eqref{fourier_first} and \eqref{second}, we obtain representation \eqref{ivanov_formula} of function $f$ as two consecutive integrals.
Due to \eqref{mestimate} we estimate Fubini's theorem applies to the consecutive integrals in the representaion \eqref{ivanov_formula} of function $f,$ and we can consider those consecutive integrals as a double integral.
\section{ Estimates for two-dimensional indicator after Ivanov.}\label{ivanov_proof}
\begin{figure}\label{fig_c}
	\centering
	\includegraphics{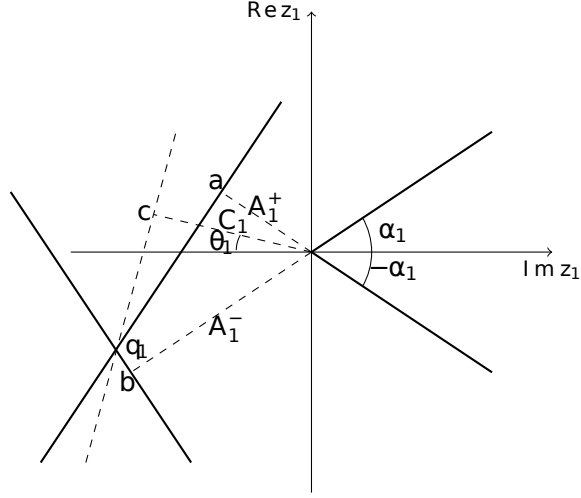}\caption{Construction of $C_1$}
\end{figure}

\subsection{ A calculation relating to figure \ref{fig_c} on page \pageref{fig_c}.}
We justify that the length on interval $|0c|$ in figure \ref{fig_c} on page \pageref{fig_c} indeed equals $C_1.$
Observing the three triangles $0qa, 0qb, 0qc$ we have
\begin{align}\label{c1}
&
Re\left(qe^{i\alpha_1}\right)=-|0a|,
\\\nonumber
&
Re\left(qe^{-i\alpha_1}\right)=-|0b|,
\\\nonumber
&
Re\left(qe^{i\theta_1}\right)=-|0c|.
\end{align}
Denote $x=Re(q),y=Im(q).$ Then we can paraphrase \eqref{c1} as
\begin{align*}&
    x\cos\left(\alpha_1\right)-y\sin\left(\alpha_1\right)=-A^+_1,
    \\\noindent &
       x\cos\left(\alpha_1\right)+y\sin\left(\alpha_1\right)=-A^-_1,
        \\\noindent &
        x\cos\left(\theta_1\right)-y\sin\left(\theta_1\right)=-|[0c]|.
\end{align*}
So that
\begin{equation*}
    |0c|=\frac{A^+_1}{2}
    \left[\frac{\sin\left(\theta_1\right)}{\sin\left(\alpha_1\right)}+\frac{\cos\left(\theta_1\right)}{\cos\left(\alpha_1\right)}\right]+\frac{A^-_1}{2}\left[\frac{\cos\left(\theta_1\right)}{\cos\left(\alpha_1\right)}-\frac{\sin\left(\theta_1\right)}{\sin\left(\alpha_1\right)}\right].
\end{equation*}
Consequently, due to \ref{c} we have $|0c|=C_1.$
\subsection{ An auxiliary estimate.}
\begin{figure}\label{fig_lambda}
	\begin{minipage}{0.45\linewidth}
		\includegraphics{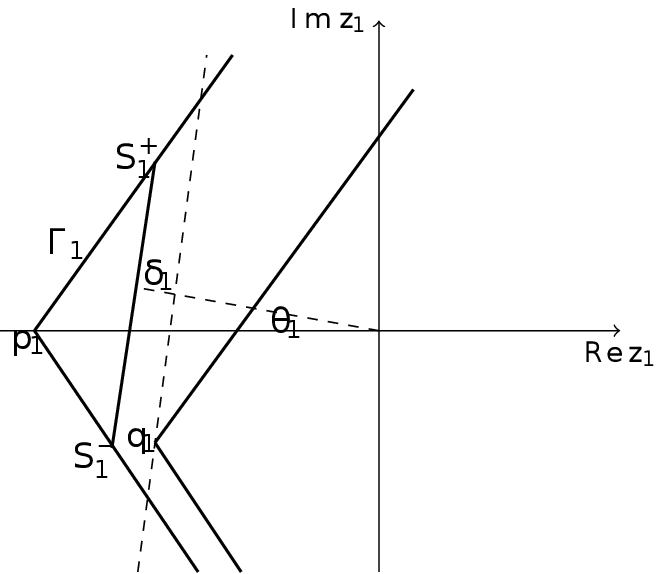}
	\end{minipage}
	\hfill
	\begin{minipage}{0.45\linewidth}
		\includegraphics{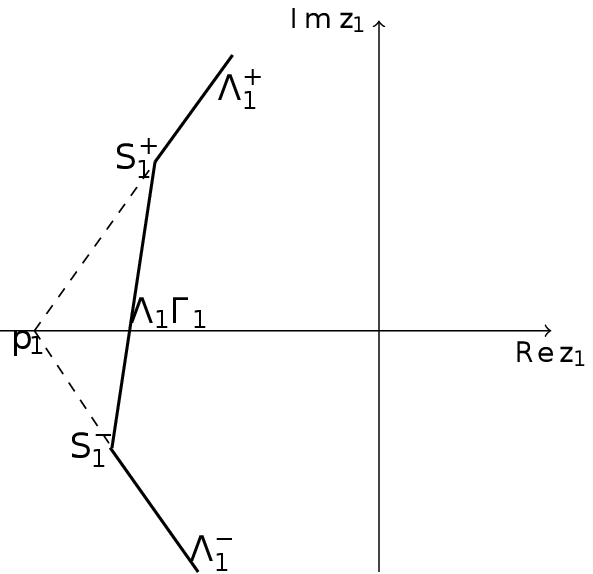}
	\end{minipage}
	\caption{Construction of $\Lambda_1$}
\end{figure}

We now estimate the integral
\begin{equation}
      \int_{\Lambda_1}e^{-Re\left(\omega_1 z_1\right)}\left|d\omega_1\right|, \quad \text{for }\arg\left(z_1\right)=\theta_1,
\end{equation}
where the curve $\Lambda_1$ is defined by figure \ref{fig_lambda} on page \pageref{fig_lambda}.
Due to the choice  \ref{c} of the constant $C_1$ and parametrization \eqref{gamma_1} of the curve $\Gamma_1,$ the curve $\Lambda_1$ is a union of three subcurves: the finite segment $\overline{\Lambda_1\setminus\Gamma_1}$ and two infinite segments $\Lambda_1^+$ and $\Lambda_1^-$ correspondingly parameterized by
\begin{align}\label{gammapart}
&
    \gamma_1(t),\quad t\geq t_+>0,
\\\nonumber &
    \gamma_1(t),\quad t\leq t_-<0,
\end{align}
where $t_+$ and $t_-$ are determined by $s_-=\gamma_1(t_-),s_+=\gamma_1(t_+).$
\paragraph{Step 1. }Note that due to the choice \ref{c} of the constant $C_1,$
\begin{equation}\label{magic}
Re\left(\omega_1 z_1\right)=Re\left(\omega_1 e^{i\theta_1}\right)\left|z_1\right|
=
-\left(C_1+\delta\right)\left|z_1\right|
,\quad\text{for }\omega_1\in\bar{\Lambda_1\setminus \Gamma_1},\;\arg\left(z_1\right)=\theta_1.
\end{equation}
And consequently
\begin{equation}\label{interval}
    \int_{\Lambda_1\setminus\Gamma_1}
    e^{-Re\left(\omega_1 z_1\right)}
    \left|d\omega\right|
    =
    e^{\left(C_1+\delta\right)\left|z_1\right|}
    \left|\Lambda_1\setminus\Gamma_1\right|.
\end{equation}
\paragraph{ Step 2. }We rewrite
\begin{equation}\label{step21}
    \int_{\Lambda_1^+}e^{-Re\left(\omega_1 z_1\right)}\left|d\omega_1\right|
    \overset{\eqref{gammapart},\eqref{gamma_1}}{=}
     \int_{t_+}^{+\infty}e^{-Re\left(\gamma_1(t)z_1\right)}dt
     \end{equation}
     Note that 
     \begin{equation*}
         \left|\theta_1\right|{<}\alpha_1{<}\pi/2.
     \end{equation*}
     Consequently
     \begin{equation*}
    Re\left(ie^{-i\alpha_1}z_1\right)>0,\quad \text{for }\arg\left(z_1\right)=\theta_1.
     \end{equation*}
     So that the following integral is convergent and equals
     \begin{equation}\label{step22}
         \int_{t_+}^{+\infty}e^{-Re\left(\gamma_1(t)z_1\right)}dt =
     \frac{e^{-Re\left(\gamma_1\left(t_+\right)z_1\right)}}{Re\left(ie^{-i\alpha_1}z_1\right)},
     \quad \text{for }\arg\left(z_1\right)=\theta_1.
\end{equation}
Note that $\gamma_1(t_+)\in\overline{\Lambda_1\setminus \Gamma_1}.$ Consequently, due to \eqref{magic}, we have 
\begin{equation}\label{step23}
   Re\left(\gamma_1\left(t_+\right)z_1\right)=-\left(C_1+\delta\right)\left|z_1\right|,
     \quad \text{for }\arg\left(z_1\right)=\theta_1.
\end{equation}
By combining \eqref{step21},\eqref{step22},\eqref{step23} we get
\begin{equation}\label{step24}
     \int_{\Lambda_{1}^+}e^{-Re\left(\omega_1 z_1\right)}\left|d\omega_1\right|\leq \frac{e^{\left(C_1+\delta\right)\left|z_1\right|}}{Re\left(ie^{-i\alpha_1}z_1\right)},   \quad \text{for }\arg\left(z_1\right)=\theta_1.
\end{equation}
Similarly,
\begin{equation}\label{step25}
    \int_{\Lambda_1^-}e^{-Re\left(\omega_1 z_1\right)}\left|d\omega_1\right|\leq \frac{e^{\left(C_1+\delta\right)\left|z_1\right|}}{Re\left(-ie^{i\alpha_1}z_1\right)},   \quad \text{for }\arg\left(z_1\right)=\theta_1.
\end{equation}
\paragraph{ Step 3.}
By combining the estimates \eqref{interval},\eqref{step24},\eqref{step25}, we get
\begin{equation}\label{step31}
     \int_{\Lambda_1}e^{-Re\left(\omega_1 z_1\right)}\left|d\omega_1\right|\leq k_{\alpha_1,\theta_1}\frac{e^{\left(C_1+\delta\right)\left|z_1\right|}}{\left|z_1\right|}, \quad \text{for }\arg\left(z_1\right)=\theta_1,
\end{equation}
where the number $k_{\alpha_1,\theta_1},$ while depending on $\alpha$ and $\theta,$ does not depend on $|z|.$
 We remark that the estimate \eqref{step31} that we have obtained is similar to the one in lemma 3.1 \cite{V}.
\subsection{ Proof of trigonometric convexity.} 
As the function $m$ is analytic on the set $\Omega_1\times\Omega_2,$ by consecutive applications of the Cauchy integral theorem in variables $z_1$ and $z_2,$
we can rewrite \eqref{ivanov_formula} as
\begin{equation}
     f\left(z_1,z_2\right)
=
\int_{\Lambda_1}\int_{\Lambda_2}
m\left(\omega_1,\omega_2\right)
e^{-\omega_1 z_1-\omega_2z_2}
d\omega_2d\omega_1,
\end{equation}
where the curve $\Lambda_1$ is constructed by figure \ref{fig_lambda}  on page \pageref{fig_lambda}.
Due to construction of the curve $\Lambda_1,$ it is bounded away from the boundary $\partial \Omega_1.$ Similarly, the curve $\Lambda_2$ is bounded away from the boundary $\partial \Omega_2.$ Consequently, due remark \ref{boundedaway}, the function $m$ is uniformly bounded (by some constant $k_{\delta}(m)$) on the Cartesian product $\Lambda_1\times\Lambda_2.$ We estimate
\begin{align}\label{nearlythere}
&
\left| f\left(z_1,z_2\right)\right|
\overset{\eqref{ivanov_formula}}{\leq}
\int_{\Lambda_1}\int_{\Lambda_2}
\left|
m\left(\omega_1,\omega_2\right)
\right|
e^{-Re(\omega_1 z_1)-Re\left(\omega_2z_2\right)}
\left|d\omega_2\right|\left|d\omega_1\right|\leq
\\\nonumber 
&
\overset{\eqref{mestimate}}{\leq}
k_{\delta}(m) \left[\int_{\Lambda_1}
e^{-Re(\omega_1 z_1)}
\left|d\omega_1\right|
\right]\cdot
\left[
\int_{\Lambda_2}
e^{-Re\left(\omega_2z_2\right)}
\left|d\omega_2\right|\right]\leq
\\\nonumber & \overset{\eqref{step31}}{\leq}
k_{\delta}(m) k_{\alpha_1,\theta_1}k_{\alpha_2,\theta_2}\frac{e^{\left(C_1+\delta\right)\left|z_1\right|}}{\left|z_1\right|}\frac{e^{\left(C_2+\delta\right)\left|z_2\right|}}{\left|z_2\right|}
,
\\\nonumber
&
\quad\text{for }    \arg\left(z_1\right)=\theta_1,
   \;
    \arg\left(z_2\right)=\theta_2.
   \end{align}
   Additionally, due to the estimate \eqref{ab} the function $f$ is bounded in: the intersection of $\Omega_1\times \Omega_2,$ and a vicinity of $0.$ By combining this fact with the estimate \eqref{nearlythere}, we obtain the following estimate on the function $f:$
   \begin{align*}
       &
       \left| f\left(z_1,z_2\right)\right|
\leq 
k_{\delta,\theta_1,\theta_2,\alpha_1,\alpha_2} e^{(C_1+\delta)\left|z_1\right|+\left(C_2+\delta\right)\left|z_2\right|},
\\\nonumber
&
\text{for }    \arg\left(z_1\right)=\theta_1,
   \;
    \arg\left(z_2\right)=\theta_2,
   \end{align*}
   where the constant $k_{\delta,\theta_1,\theta_2,\alpha_1,\alpha_2},$ 
   while depending on $\delta,\theta_1,\theta_2,\alpha_1,\alpha_2,$
   does not depend on $\left|z_1\right|$ or $\left|z_2\right|.$
\section{ Proof of sharpness.}\label{sharp_section}
\begin{figure}\label{fig_sharp}
\centering
\includegraphics{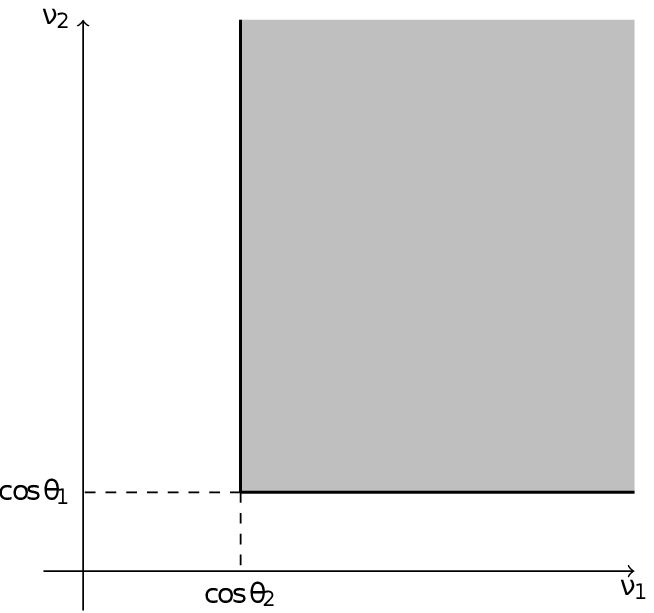}
\caption{$T_f\left(\theta_1,\theta_2\right)$}
\end{figure}

Consider the entire function
 $f(z_1,z_2)=e^{z_1+z_2}.$ The function $f$ is of exponential type,
 \begin{equation*}
 \left|f\left(z_1,z_2\right)\right|\leq e^{\left|z_1\right|\cos \theta_1+\left|z_2\right|\cos \theta_2},\quad\text{for }z_1,z_2\in\mathbb{C}.    
 \end{equation*}
For our choice of the function $f,$ the set $T_f(\theta_1, \theta_2)$ defined by \eqref{T} equals
\begin{align}\label{ivanov_evaluated}
T_f\left(\theta_1, \theta_2\right)
=
\left\lbrace(\nu_1,\nu_2) : \nu_1\geq\cos \theta_1, \nu_2\geq\cos \theta_2   \right\rbrace
\end{align}
(see figure \ref{fig_sharp} on page \pageref{fig_sharp}).
In particular, for $\alpha_1=\alpha_2=\pm\frac{\pi}{4}$ we have
\begin{equation*}
T_f\left(\pm\frac{\pi}{4}, \pm\frac{\pi}{4}\right)
=
\left\lbrace(\nu_1,\nu_2) \colon \nu_1\geq \frac{\sqrt{2}}{2}, \nu_2\geq\frac{\sqrt{2}}{2}\right\rbrace.
\end{equation*}
Thus
\begin{equation*}
A^+_1=A^+_2=A^-_1=A^-_2=\frac{\sqrt{2}}{2}. 
\end{equation*}
Take $\theta_1=\theta_2 = 0.$ By \eqref{c} we evaluate
\begin{equation*} 
C_1=C_2=\frac{\sqrt{2}}{4} \frac{2}{\sqrt{2}}+\frac{\sqrt{2}}{4} \frac{2}{\sqrt{2}} = 1.
\end{equation*}
Thus, according to theorem \ref{ivanov} the following inequality holds:
\begin{equation*}
f(z_1, z_2) \leq k_{\varepsilon} e^{(1+\varepsilon)\left|z_1\right|+(1+\varepsilon)\left|z_2\right|},\quad\text{for }\arg z_1 = \arg z_2 = 0.
\end{equation*}
On the other hand, by \eqref{ivanov_evaluated}
\begin{equation*} 
T_f(0, 0)=\left\lbrace(\nu_1,\nu_2) : \nu_1\geq 1, \nu_2\geq 1  \right\rbrace.
\end{equation*}
That is, in this case, theorem \ref{ivanov} is sharp.
\section*{Acknowledgments}
The first author  was  funded  by a grant of the Russian Science Foundation (project No. 20-11-20117)

The second author was partially funded by a project of the Russian Ministry of Science and Education for the organization and development of scientific and educational mathematical centers (agreement No. 075-02-2022-893). 

\bibliography{tcmii}

\section*{}
\begin{minipage}{0.5\linewidth}
    Aleksandr Mkrtchyan\\
AMkrtchyan@sfu-kras.ru\\
    24/5 Marshal Bagramian ave.\\
    Yerevan, 0019, Republic of Armenia
\end{minipage}
\hfill
\begin{minipage}{0.5\linewidth}
    Armen Vagharshakyan\\
    avaghars@kent.edu \\
24/5 Marshal Bagramian ave.\\
    Yerevan, 0019, Republic of Armenia
\end{minipage}

\end{document}